\documentclass[11pt,reqno]{amsart}
\pdfoutput=1

\usepackage{amssymb}

\usepackage{graphicx}
\usepackage{enumerate}

\catcode`\@=11

\long\def\@savemarbox#1#2{\global\setbox#1\vtop{\hsize\marginparwidth 
  \@parboxrestore\tiny\raggedright #2}}
\newcommand\lref[1]{\ref{#1}%
\@ifundefined{r@DisplaY #1}{}{ (#1)}}

\newcommand\fakelabel[2]{\@bsphack\if@filesw {\let\thepage\relax
   \newcommand\protect{\noexpand\noexpand\noexpand}%
\xdef\@gtempa{\write\@auxout{\string
      \newlabel{#1}{{#2}{\thepage}}}}}\@gtempa
   \if@nobreak \ifvmode\nobreak\fi\fi\fi\@esphack}

\def\SL@margintext#1{{\showlabelsetlabel{\tiny\{\SL@prlabelname{#1}\}}}}
\catcode`\@=12

\def\Empty{}
\newcommand\oplabel[1]{
  \def\OpArg{#1} \ifx \OpArg\Empty {} \else
        \label{#1}
  \fi}
\newtheorem{theoremSt}{Theorem}[section]

\newtheorem{exampleSt}[theoremSt]{Example}
\newtheorem{exerciseSt}[theoremSt]{Exercise}

\newcommand\MakeStEnv[1]{
  \newenvironment{#1}[1]{
  \begin{#1St} \oplabel{##1}%
  \global\def\CrntSt{\thetheoremSt}%
}{ 
  \end{#1St} }
  \newenvironment{#1+}[1]{
  \begin{#1St} \label{##1}%
  \label{DisplaY ##1}%
  \global\def\CrntSt{\thetheoremSt}%
  \def\Labl{##1}\ifx\Labl\Empty{} \else {\em (\Labl)\,}\fi%
}{ 
  \end{#1St} }
}
\MakeStEnv{theorem}
\MakeStEnv{corollary}
\MakeStEnv{proposition}
\MakeStEnv{lemma}
\MakeStEnv{definition}
\MakeStEnv{conjecture}
\MakeStEnv{problem}
\MakeStEnv{question}

\newlength{\saveu}

\newenvironment{pf*}[1]{%
 \begin{proof}[#1]%
}{ 
 \end{proof}
}

\newcommand{\finishproof}[1]{ 
  \def\FPArg{#1}
  \ifx\FPArg\Empty
        \newcommand\FPArg{\CrntSt}  \fi
  \smallbreak\noindent\makebox[\textwidth]{\hfill\fbox{\FPArg}}
  \medbreak\noindent
}

\newcommand\FF{{\mathcal F}}

\newcommand\LL{{\mathcal L}}
\newcommand\MM{{\mathcal M}}

\newcommand\PP{{\mathcal P}}

\newcommand\PMF{{\PP\kern-2pt\MM\FF}}

\newcommand\PML{{\PP\kern-2pt\MM\LL}}

\newcommand\bbR{{\mathord{\text{I\kern-2pt R}}}}        
\newcommand\bbH{{\mathord{\text{I\kern-2pt H}}}}        

\newcommand\bigrightarrow[1]{\hbox to #1{\rightarrowfill}}
\newcommand\bigleftarrow[1]{\hbox to #1{\leftarrowfill}}

\newcommand\semidir{\mathrel{\hbox{\vrule depth-.03ex height1.1ex\kern-0.15em$\times$}}}

\numberwithin{equation}{section}

\begin{document}

\title{Dynamics on character varieties: a survey}

\author{Richard D. Canary}
\address{University of Michigan}
\date{\today}
\thanks{The author gratefully acknowledge support from the National
  Science Foundation}

\begin{abstract}
We survey recent work on the dynamics of  the outer automorphism group ${\rm Out}(\Gamma)$
of a word hyperbolic group on spaces of (conjugacy classes of) representations of $\Gamma$ into
a semi-simple Lie group $G$. All these results are motivated by the fact that the mapping class
group  of a closed surface acts properly discontinuously on the Teichm\"uller space of the surface.

We explore two settings. 1) The work of Guichard, Labourie and Wienhard establishing
the proper discontinuity of the action of outer automorphism groups on  spaces of Anosov
representations. 2) The work of Canary, Gelander, Lee, Magid, Minsky and Storm on the case where
$\Gamma$ is the fundamental group of a compact 3-manifold with boundary and $G={\rm PSL}(2,\mathbb C)$.
\end{abstract}

\maketitle

\section{Overview}

It is a classical result, due to Fricke, that the mapping class group ${\rm Mod}(F)$ of a closed oriented surface $F$ acts properly discontinuously
on the Teichm\"uller space $\mathcal{T}(F)$ of marked  hyperbolic (or conformal) structures on $F$. 
The mapping class group ${\rm Mod}(F)$ may be identified with an index two subgroup of
the outer automorphism group ${\rm Out}(\pi_1(F))$ of $\pi_1(F)$ and Teichm\"uller space may
be identified with a component of the space $X(\pi_1(F),{\rm PSL}(2,\mathbb R))$ of conjugacy classes
of representations of $\pi_1(F)$ into ${\rm PSL}(2,\mathbb R)$.

It is natural to attempt to generalize this beautiful result  by studying the dynamics of the action of
the outer automorphism group 
${\rm Out}(\Gamma)$ of a finitely generated group $\Gamma$ on the space $X(\Gamma,G)$
of (conjugacy classes of) representations
into a semi-simple Lie group $G$. We recall that 
$${\rm Out}(\Gamma)={\rm Aut}(\Gamma)/{\rm Inn}(\Gamma)$$
where ${\rm Inn}(\Gamma)$ is
the group of inner automorphisms of $\Gamma$. We let
$$X(\Gamma,G)={\rm Hom}(\Gamma,G)/G.$$
In the case that $G$ is a complex semi-simple Lie group, e.g. ${\rm PSL}(2,\mathbb C)$,
one may instead take the Mumford
quotient (see, for example, Lubotzky-Magid \cite{lubotzky-magid}),  which  has the structure of an algebraic variety.
This additional structure will not be relevant to our considerations. However, we will abuse notation in
the traditional manner and refer to $X(\Gamma,G)$ as the $G$-character variety of $\Gamma$.

The first natural generalization is that if $\Gamma$ is word hyperbolic, then ${\rm Out}(\Gamma)$ acts properly
discontinuously on the space ${\rm CC}(\Gamma,{\rm PSL}(2,\mathbb C))$ of convex cocompact
representations of $\Gamma$ into ${\rm PSL}(2,\mathbb C)$. We recall that 
${\rm CC}(\Gamma,{\rm PSL}(2,\mathbb C))$ is an open subset of $X(\Gamma,{\rm PSL}(2,\mathbb C))$,
so this shows that ${\rm CC}(\Gamma,{\rm PSL}(2,\mathbb C))$ is a domain of discontinuity for the
action of ${\rm Out}(\Gamma)$. More generally, if $G$ is a rank one semi-simple Lie group, then
${\rm Out}(\Gamma)$ acts properly discontinuously on the subset $CC(\Gamma, G)$ of convex cocompact
representations of $\Gamma$ into $G$ (see Section \ref{convex cocompact}).

Labourie \cite{labourie-anosov} introduced the notion of an Anosov representation into
a semi-simple Lie group $G$  and the theory was further developed by Guichard and Wienhard
\cite{guichard-wienhard-domains}.
One may view Anosov representations as
the natural analogue of convex cocompact representations in the setting of higher rank Lie groups. In fact,
if $G$ has rank one, then a representation is Anosov if and only if it is convex cocompact.
(Kleiner-Leeb \cite{kleiner-leeb} and Quint \cite{quint} showed that there are no ``interesting'' convex 
cocompact representations into higher rank Lie groups.)
The space of Anosov representations is an open subset of $X(\Gamma,G)$ and Labourie \cite{labourie-anosov}
and Guichard-Wienhard \cite{guichard-wienhard-domains} showed that ${\rm Out}(\Gamma)$ acts
propertly discontinuously on the space of Anosov representations. 

Finally,  we discuss work of  Canary, Lee, Magid, Minsky and Storm
(\cite{canary-magid,canary-storm1,canary-storm2,lee-Ibundle,lee-compbody,minsky-PS})
which exhibits domains of discontinuity for
the action of ${\rm Out}(\Gamma)$ on $X(\Gamma,{\rm PSL}(2,\mathbb C))$ which are (typically)
strictly larger than ${\rm CC}(\Gamma,{\rm PSL}(2,\mathbb C))$. Often, these domains of discontinuity include
representations which are neither discrete nor faithful.
In the case of the free group $F_n$ (with $n\ge 3$), Gelander and Minsky \cite{gelander-minsky},
have found  an open subset of $X(F_n,{\rm PSL}(2,\mathbb C))$ where ${\rm Out}(\Gamma)$
acts ergodically.

We have restricted our discussion to the case where $G$ has no compact factors, but we note that when
$G$ is compact and $\Gamma$ is either a free group or the fundamental group of a closed  hyperbolic
surface, then the action of ${\rm Out}(\Gamma)$ on $X(\Gamma,G)$ is known to be ergodic in many cases,
see, for example, Goldman \cite{goldman-su2,goldman-free}, Pickrell-Xia \cite{pickrell-xia}, Palesi \cite{palesi}
and Gelander \cite{gelander-free}.
For a discussion of representations of surface groups with a somewhat different viewpoint, which also
includes a discussion of the case when $G$ is compact, we recommend Goldman \cite{goldman-survey}.

One would ideally like a dynamical dichotomy for the action of ${\rm Out}(\Gamma)$ on $X(\Gamma,G)$, i.e.
a decomposition into an open set where the action is properly discontinuous and a closed set where
the action is chaotic in some precise sense.  The only case where $G$ is non-compact and
such a dynamical dichotomy is understood
is the case where $\Gamma=F_2$ is the free group on 2 generators and $G={\rm SL}(2,\mathbb R)$ (see
Goldman \cite{goldman-pt}). We discuss this work briefly at the end of section \ref{teich}.

\medskip\noindent
{\bf Acknowledgements:} I would like to thank Michelle Lee for her thoughtful comments on an early
version of the paper. I would also like to thank Mark Hagen for several useful conversations concerning
the proof of Proposition \ref{gens-suffice}.

\section{Basic definitions}

Throughout this article $G$ will be a semi-simple Lie group with trivial center and no compact
factors and $\Gamma$ will be a torsion-free word hyperbolic group which is not virtually cyclic.
If $K$ is a maximal compact subgroup of $G$, then $X=G/K$ is a symmetric space
and $G$ acts as a group of isometries of $X$.  The symmetric space is non-positively
curved and is negatively curved if and only if $G$ has rank one.
The most basic example is when $G={\rm PSL}(2,\mathbb R)$
and $X=\mathbb H^2$. Similarly, if $G={\rm PSL}(2,\mathbb C)$, then $X=\mathbb H^3$, and more generally,
if $G=SO_0(n,1)$, then $X=\mathbb H^n$. All these examples have rank one. The simplest example
of a higher rank semi-simple Lie group is ${\rm PSL}(n,\mathbb R)$ (when $n\ge 3$)
which has rank $n-1$. For the purposes of this article, it will suffice to restrict one's attention
to these examples.

We consider the subset $AH(\Gamma,G)\subset X(\Gamma,G)$ of discrete, faithful representations.
If $\rho\in AH(\Gamma,G)$, then $N_\rho=X/\rho(G)$ is a manifold which is a locally symmetric space modeled on $X$
with fundamental group isomorphic to $\Gamma$. One may then think of $AH(\Gamma,G)$ as the space
of marked locally symmetric spaces modelled on $X$ with fundamental group isomorphic to $\Gamma$. So,
$AH(\Gamma,G)$ is one natural generalization of the Teichm\"uller space of a closed oriented surface $F$, which
is the space of marked hyperbolic structures on $F$. It is a  consequence of the Margulis
Lemma that $AH(\Gamma,G)$ is  a closed subset of $X(\Gamma,G)$ (see Kapovich \cite[Thm. 8.4]{kapovich}).

If $\rho:\Gamma\to G$ is a representation and $x\in X$, then we can define the 
{\em orbit map} $\tau_{\rho,x}:\Gamma\to X$ defined by $\tau_{\rho,x}(g)=\rho(g(x))$.
The orbit map is said to be a {\em quasi-isometric embedding} if there
exists $K$ and $C$ such that
$$\frac{d(g,h)}{K}-C\le d(\tau_{\rho,x}(g),\tau_{\rho,x}(h))\le K d(g,h)+C$$
for all $g,h\in \Gamma$, where $d(g,h)$ is the word length of $gh^{-1}$ with respect to some
fixed finite  generating set of $\Gamma$.  It is easily checked that  if the orbit map
is a quasi-isometric embedding for a single fixed point $x\in X$ and choice of generating set for $\Gamma$,
then the orbit map associated to any other point is also a quasi-isometric embedding with
respect to any finite generating set for $\Gamma$.

After one has chosen a finite generating set for $\Gamma$, one may form the associated
Cayley graph $C_\Gamma$. The orbit map then extends to a map
\hbox{$\bar\tau_{\rho,x}:C_\Gamma\to X$} by simply taking any edge to the geodesic joining the images of
its endpoints.  The extended orbit map is a quasi-isometric embedding if and only if the original
orbit map was a quasi-isometric embedding.

If $\gamma\in\Gamma$, then we let $\ell_\rho(\Gamma)$ denote the translation length
of $\rho(\gamma)$, i.e.
$$\ell_\rho(\gamma)=\inf_{x\in X} d_X(\rho(\gamma)(x),x).$$ 
Let $||\gamma||$ denote the translation length of the action of $\gamma$ on $\Gamma$
(or on the Cayley graph  $C_\Gamma$) with respect to our fixed generating set for $\Gamma$.
Equivalently, $||\gamma||$ is the minimal word length of an element conjugate to $\gamma$.
We say that the representation $\rho:\Gamma\to G$ is {\em well-displacing} if there exists $J$ and $B$ such that
$$\frac{1}{J}||\gamma||-B\le \ell_\rho(\gamma)\le J||\gamma||+B$$
for all $\gamma\in\Gamma$.

We observe that if $G$ has rank one, then $\rho$ is well-displacing if $\bar\tau_{\rho,x}$ is a quasi-isometric embedding.
This fact is used in the proofs of Fricke's Theorem and its generalization, Theorem \ref{rank one}.
It also plays a role in Minsky and Lee's results, described in section \ref{freepsl2c}.

\begin{proposition}{qigiveswd}{}
Suppose that $\Gamma$ is a word hyperbolic group and $G$ is a simple rank one Lie group with
trivial center and no compact factors.
If $[\rho]\in X(\Gamma,G)$ and  the extended orbit
map $\bar\tau_{\rho,x}$ is a quasi-isometric embedding for some (any) $x\in X$, then 
$\rho$ is well-displacing.
\end{proposition}

\noindent{\em Sketch of proof:}
For simplicity, we will  assume that every $\gamma\in\Gamma$ has an axis in $C_\Gamma$, i.e. a geodesic $L_\gamma$ in
$C_\Gamma$ such that $\gamma(L_\gamma)=L_\gamma$ and if $y\in L_\gamma$, then $d(y,\gamma(y))=||\gamma||$.
Suppose that $\bar\tau_{\rho,x}$ is a $(K,C)$-quasi-isometric embedding. Since $G$ has rank one, $X$ has sectional
curvature bounded above by $-1$. Suppose that $\gamma\in\Gamma$ and $L_\gamma$ is an axis for $\gamma$ on $C_\Gamma$.
Then, the restriction of $\tau_{\rho,x}$ to $L_\gamma$ is a $(K,C)$-quasi-isometric embedding, i.e.
$\tau_{\rho,x}|_{L_\gamma}$ is a $(K,C)$-quasi-geodesic. The fellow traveller property for $CAT(-1)$-spaces (see for example,
Bridson-Haefliger \cite[Thm III.H.1.7]{bridson-haefliger}) implies
that there exists $D=D(K,C)$ such that  the Hausdorff distance between $\bar\tau_{\rho,x}(L_\gamma)$ and the axis $A_\gamma$ for the action of
$\rho(\gamma)$ on $X$ is at most $D$. If we pick $z\in A_\gamma$ and $y\in L_\gamma$ so that $d(z,\bar\tau_{\rho,x}(y))\le D$,
then, since $\bar\tau_{\rho,x}$ is $\rho$-equivariant
$$|d(\bar\tau_{\rho,x}(y),\bar\tau_{\rho,x}(\gamma(y)))-d(z,\rho(\gamma(z)))|\le 2D.$$
Since $\bar\tau_{\rho,x}$ is a $(K,C)$-quasi-isometric embedding, it follows that 
$$\frac{1}{K}d(y,\gamma(y))-(C+2D)\le d(z,\rho(\gamma)(z))\le Kd(y,\gamma(y))+(C+2D).$$
Since, $||\gamma||=d(y,\gamma(y))$ and $\ell_\rho(\gamma)=d(z,\rho(\gamma)(z))$, it follows that
$\rho$ is $(K,C+2D)$-well-displacing.

In general, not every element of $\gamma$ need have an axis, but there always exists $(L,A)$ so that
every element has a $(L,A)$-quasi-axis, i.e. a $(L,A)$-quasi-geodesic $L_\gamma$ in $C_\Gamma$ such that
if $y\in L_\gamma$, then $d(y,\gamma))=||\gamma||$ (see Lee \cite{lee-thesis} for a proof of the existence of
quasi-axes).  Once one has made this observation
the argument in the previous paragraph immediately generalizes. This completes the sketch of the proof.

\medskip

More generally, Delzant, Guichard, Labourie and Mozes \cite{labourie-displacing} proved that in our
context, a representation is well-displacing if and only if its orbit maps are quasi-isometries.

\begin{theorem}{wdandqi}{\em (Delzant-Guichard-Labourie-Mozes \cite{labourie-displacing})}
If $\Gamma$ is a word hyperbolic group and $G$ is a semi-simple Lie group with
trivial center and no compact factors, then $\rho\in X(\Gamma,G)$ is well-displacing if and only the orbit
map $\tau_{\rho,x}$ is a quasi-isometric embedding for some (any) $x\in X$. 
\end{theorem}

We will need the following observation about outer automorphism groups of word hyperbolic
groups, which will be established in an appendix to the paper.

\begin{proposition}{gens-suffice}
If $\Gamma$ is a torsion-free word hyperbolic group, then there exists a finite set $\mathcal B$ of elements
of $\Gamma$ such that for any $K$,
$$\{\phi\in {\rm Out}(\Gamma)\ |\  ||\phi(\beta)||\le K \ \ {\rm if}\ \beta\in{\mathcal{B}} \}$$
is finite.

Moreover, if $\Gamma$ admits a convex cocompact action on the symmetric space associated
to a rank one Lie group and has finite generating set $S$ one may take
$\mathcal B$ to consist of all the elements of $S$ and all products of two elements of $S$.
\end{proposition}

\noindent
{\bf Remark:} Proposition \ref{gens-suffice} is a generalization of earlier results by Minsky
\cite{minsky-PS}, for free groups, and Lee \cite{lee-compbody} for free products of
infinite cyclic groups and closed orientable hyperbolic surface groups. See also
Bogopolski-Ventura \cite{bogo-vent} for a related result. It seems likely that Proposition
\ref{gens-suffice} can also be proven using their techniques.

\section{Teichm\"uller space}
\label{teich}

In this section, we sketch a proof of Fricke's Theorem that the mapping class group  of a closed oriented
surface acts
properly discontinuously on its Teichm\"uller space. This basic proof will
be adapted to a variety of more general settings.

If $\rho\in AH(\pi_1(F),{\rm PSL}(2,\mathbb R))$, then $\rho$ is discrete and faithful
and the quotient $N_\rho=\mathbb H^2/\rho(\pi_1(F))$ is a hyperbolic
surface homotopy equivalent to $F$. Moreover, there exists a homotopy equivalence $h_\rho:F\to N_\rho$
in the homotopy class determined by $\rho$. The Baer-Nielsen theorem implies that $h_\rho$ is homotopic
to a homeomorphism.

If $\rho\in AH(\pi_1(F),{\rm PSL}(2,\mathbb R))$, then, since $N_\rho$ is closed, the Milnor-\v{S}varc Lemma 
(see, for example, Bridson-Haefliger \cite[Prop. I.8.19]{bridson-haefliger}) implies
that any orbit map $\tau_{\rho,x}:\pi_1(F)\to \mathbb H^2$ is  a quasi-isometric embedding.
On the other hand, if any orbit map is a quasi-isometric embedding, then $\rho$ is discrete and
faithful. So, $\rho\in AH(\pi_1(F),{\rm PSL}(2,\mathbb R))$ if and only if any orbit map is a quasi-isometric embedding.

One can check that if $\tau_{\rho,x}$
is a quasi-isometric embedding, then there exists an open neighborhood  $U$ of $\rho$ so that
if $\rho'\in U$, then the orbit map $\tau_{\rho',x}$ is also a quasi-isometric embedding,
so $AH(\pi_1(F),{\rm PSL}(2,\mathbb R))$ is open. We will outline a proof of a more general statement later
(see Proposition \ref{openinrank1}), but for the
moment we note that this is a rather standard consequence of the stability of quasi-geodesics in
negatively curved spaces.  Since we have already observed that $AH(\pi_1(F),{\rm PSL}(2,\mathbb R))$ is closed,
it is a collection of components of $X(\pi_1(F),{\rm PSL}(2,\mathbb R))$.

In fact, $AH(\pi_1(F),{\rm PSL}(2,\mathbb R))$ has two components, one of which is identified 
with ${\mathcal{T}}(F)$ and the other
of which is identified with ${\mathcal{T}}(\bar F)$ where $\bar F$ is $F$ with the opposite orientation
(see Goldman \cite{goldman-components}).
(Here $\rho\in{\mathcal{T}}(F)$ if and only if $h_\rho$ is orientation-preserving.)
Each component is homeomorphic to $\mathbb R^{6g-6}$ where $g$ is the genus of $F$.
${\rm Out}(\pi_1(F))$ may be identified with the space of (isotopy classes of) homeomorphisms of $F$ 
and the mapping class group ${\rm Mod}(F)$ is simply the  index two
subgroup consisting of orientation-preserving homeomorphisms. 
(See Farb-Margalit  \cite{farb-margalit} for a discussion of the mapping class group and its action on Teichm\"uller space.)
The mapping class group
preserves each component of $AH(\pi_1(F),{\rm PSL}(2,\mathbb R))$.

\begin{theorem}{}{\em (Fricke)}  
If $F$ is a closed, oriented surface of genus $g\ge 2$, then
${\rm Out}(\pi_1(F))$ acts properly discontinuously on $AH(\pi_1(F),{\rm PSL}(2,\mathbb R))$. In particular,
${\rm Mod}(F)$ acts properly discontinuously on ${\mathcal{T}}(F)$.
\end{theorem}

\noindent
{\em Outline of proof:} 
If $\rho\in AH(\pi_1(F),{\rm PSL}(2,\mathbb R))$, then, since the orbit map is a quasi-isometric embedding,
Proposition \ref{qigiveswd} implies that $\rho$ is well-displacing, i.e. that
$\ell_\rho(\gamma)$ is $(J,B)$-comparable to $||\gamma||$ for all $\gamma\in\pi_1(F)$ and some $(J,B)$.

If $\{\phi_n\}$ is a sequence of distinct elements of ${\rm Out}(\pi_1(F))$, then
Proposition \ref{gens-suffice} implies that, up to subsequence, there
exists $\beta\in\pi_1(F)$ so that $||\phi_n^{-1}(\beta)||\to \infty$. 
Therefore,  since $\rho$ is well-displacing, 
\hbox{$\ell_\rho(\phi_n^{-1}(\beta))=\ell_{\rho\circ\phi_n^{-1}}(\beta)\to\infty$}, which implies that 
\hbox{$\phi_n(\rho)=\rho\circ\phi_n^{-1}\to\infty$} in $X(\pi_1(F),{\rm PSL}(2,\mathbb R))$.
Therefore,  the action of ${\rm Out}(\pi_1(F))$  on $AH(\pi_1(F), {\rm PSL}(2,\mathbb R))$ has discrete orbits.

In order to extend this proof to show that the action is actually properly discontinuous, it only remains
to observe that on any compact subset $R$ of $AH(\pi_1(F),{\rm PSL}(2,\mathbb C))$, one may choose $J$ and $B$
so that 
$$\frac{1}{J}||\gamma||-B\le \ell_\rho(\gamma)\le J||\gamma||+B$$
for all $\gamma\in\pi_1(F)$ and all  $\rho\in R$. 
(In order to prove this, note that Proposition \ref{openinrank1} implies that there exists
$(K,C)$ and $x\in \mathbb H^2$ so that each $\rho\in R$ has a representative such that
$\tau_{\rho,x}$ is a $(K,C)$-quasi-isometric embedding, while the proof of Proposition \ref{qigiveswd}
shows that one can then find $(J,B)$ so that each $\rho\in R$ is $(J,B)$-well-displacing.) One then uses
Proposition \ref{gens-suffice} to show that if $\{\phi_n\}$ is a sequence of distinct elements of ${\rm Out}(\pi_1(F))$,
then $\{\phi_n(R)\}$ exists every compact subset of $X(\pi_1(F),{\rm PSL}(2,\mathbb R))$.
This completes our sketch of the argument.

\medskip

The dynamics of the action of ${\rm Out}(\pi_1(F))$ on the remainder of 
$X(\pi_1(F),{\rm PSL}(2,\mathbb R))$ is more mysterious.
Goldman \cite{goldman-components} proved that $X(\pi_1(F),{\rm PSL}(2,\mathbb R))$ has $4g-3$ components
(which are enumerated by the euler number of the representations in the component). He also
made the following conjecture:

\medskip\noindent
{\bf Conjecture:} (Goldman) {\em If $F$ is a closed oriented surface of genus $g\ge 2$, then
${\rm Mod}(F)$ acts ergodically on each of  the $4g-5$ components of
$X(\pi_1(F),{\rm PSL}_2(\mathbb{R}))$ which do not correspond to either $\mathcal{T}(F)$ or $\mathcal{T}(\bar{F})$.}

\medskip
Evidence for Goldman's conjecture is provided by his analysis
 \cite{goldman-pt} of the dynamics of the action of ${\rm Out}(F_2)$ on $X(F_2,{\rm SL}(2,\mathbb R))$,
where $F_2$ is the free group on two generators, $a$ and $b$.
In this case, the action of ${\rm Out}(F_2)$ preserves the level set of the function $\kappa$ where $\kappa(\rho)$ is
the trace of $\rho([a,b])$.  He shows that for each $t\in\mathbb R$,  the level
set $\kappa^{-1}(t)$ contains an open set (possibly empty) consisting of representations associated to (possibly singular)
hyperbolic surfaces where the action is properly discontinuous and the action on the complement (again possibly empty)
is ergodic.

\section{Convex cocompact representations}
\label{convex cocompact}

In this section, we will assume that $G$ is a rank one semi-simple Lie group, e.g. ${\rm PSL}(2,\mathbb C)$, 
so that the associated symmetric space is negatively curved. In fact we may assume that its sectional
curvature is bounded above by $-1$.
If $\rho:\Gamma\to G$ is a discrete faithful representation, then $N_\rho=X/\rho(\Gamma)$
is a manifold with fundamental group $\Gamma$. 
We will say that a discrete, faithful representation  $\rho$ is {\em convex cocompact}
if $N_\rho$ contains a compact, convex submanifold whose inclusion into $N_\rho$ is a homotopy
equivalence. The Milnor-\v{S}varc Lemma again assures
us that if $\rho$ is convex cocompact, then the orbit map $\tau_{\rho,x}:\Gamma\to X$ is a quasi-isometric embedding. 
On the other hand,  if the orbit map is a quasi-isometric embedding, then
the convex hull of $\tau_{\rho,x}(\Gamma)$ lies in a bounded neighborhood of $\tau_{\rho,x}(\Gamma)$,
so has compact quotient in $N_\rho$ (see, for example, Bourdon \cite[Cor. 1.8.4]{bourdon}).
It follows that if $G$ has rank one, then $\rho$ is
convex cocompact if and only if its associated orbit maps are quasi-isometric embeddings.

Let $CC(\Gamma,G)\subset X(\Gamma,G)$ denote the set of (conjugacy classes of)
convex cocompact representation of $\Gamma$ into $G$.
We observe that $CC(\Gamma,G)$ is an open subset of $X(\Gamma,G)$. However, it is no longer
necessarily closed, so it need not be a collection of components of $X(\Gamma,G)$.

\begin{proposition}{openinrank1}
If $G$ is a rank one semisimple Lie group and $\Gamma$ is a torsion-free word hyperbolic
group, then $CC(\Gamma,G)$ is open in $X(\Gamma,G)$. 

Moreover, if $[\rho_0]\in CC(\Gamma,G)$, then
there exists an open neighborhood $U$ of $[\rho_0]$ in $X(\Gamma,G)$ and $x\in X$,  $K$ and $C$ so
that if $[\rho]\in U$, then $[\rho]$ has a representative $\rho$ such that
$\bar\tau_{\rho,x}$ is a $(K,C)$-quasi-isometric embedding.
\end{proposition}

\noindent
{\em Sketch of proof:} Suppose that $[\rho_0]\in CC(\Gamma,G)$, so, given $x\in X$, there exists $K_0$ and $C_0$ so that
$\bar\tau_{\rho_0,x}$ is a $(K_0,C_0)$-quasi-isometric embedding. We recall that, since $X$ is $CAT(-1)$,
there exists $A>0$, $K$ and $C$ so that if
$J$ is an interval in $\mathbb R$ and $h:J\to X$ is a $(K_0,3C_0)$-quasi-isometric embedding on every subsegment of $J$ of
length  at most $A$, then $h$ is a $(K,C)$-quasi-isometric embedding (see, for example, Ghys-de la Harpe \cite[Chapter 5]{ghys-delaharpe}).

Choose a neighborhood $U$ of $[\rho_0]$ so that if $[\rho]\in U$, then
it has a representative $\rho$ such that  \hbox{$d(\tau_{\rho,x}(z),\tau_{\rho_0,x}(z))\le C_0$} if $d(z,id)\le A+1$.
It follows that  the restriction of $\tau_{\rho,x}$ to every geodesic segment of length at most $A+1$
originating at the origin is a $(K_0,3C_0)$-quasi-isometric embedding. Since $\bar\tau_{\rho,x}$ is $\rho$-equivariant, it follows
that the restriction of  $\bar\tau_{\rho,x}$ to any segment in $C_\Gamma$ of length at most $A$ is a $(K_0,3C_0)$-quasi-isometric
embedding. Therefore, the restriction of $\bar\tau_{\rho,x}$ to any geodesic in $C_\Gamma$ is a $(K,C)$-quasi-isometric
embedding, which implies that $\bar\tau_{\rho,x}$ is itself a $(K,C)$-quasi-isometric
embedding.

\medskip\noindent
{\bf Remarks:} (1) We notice that the fact that $X$ is negatively curved is crucial here.
To see what can go wrong in the non-positively curved case, we recall that any translation in
Euclidean space arises as a limit of rotations (of higher and higher order).

(2) Marden \cite{marden} proved Proposition \ref{openinrank1} in the case that $G={\rm PSL}(2,\mathbb C)$,
while Thurston \cite{thurston-notes} established it in the case that $G=SO_0(n,1)$.

\medskip

One may then easily generalize our proof of Fricke's Theorem to obtain:

\begin{theorem}{rank one}
If $G$ is a rank one semisimple Lie group and $\Gamma$ is a torsion-free word hyperbolic
group, then ${\rm Out}(\Gamma)$ acts properly discontinuously on $CC(\Gamma,G)$.
\end{theorem}

If $F$ is a closed oriented surface of genus at least two, then
$CC(\pi_1(F),{\rm PSL}(2,\mathbb C))$ is the well-studied space of quasifuchsian representations.
Bers \cite{bers-su} proved that $CC(\pi_1(F),{\rm PSL}(2,\mathbb C))$ is identified with
${\mathcal{T}}(F)\times {\mathcal{T}}(\bar F)$. Brock, Canary and Minsky \cite{ELC2} proved that the
closure of $CC(\pi_1(F),{\rm PSL}(2,\mathbb C))$ is the set  $AH(\pi_1(F),{\rm PSL}(2,\mathbb C))$
of discrete, faithful representations of  $\pi_1(F)$ into ${\rm PSL}(2,\mathbb C)$, 
while Bromberg and Magid \cite{bromberg-PT,magid-lc}
proved that $AH(\pi_1(F),{\rm PSL}(2,\mathbb C))$ fails to be locally connected.

Thurston \cite{thurstonII} (see also Otal \cite{otal-fibre}) proved that if $\phi\in {\rm Mod}(F)$ is a pseudo-Anosov mapping
class, then it has a fixed point in the closure of $CC(\pi_1(F),{\rm PSL}(2,\mathbb C))$. In particular,
${\rm Mod}(F)$ does not act properly discontinuously on $AH(\pi_1(F),{\rm PSL}(2,\mathbb C))$.
Canary and Storm \cite{canary-storm1} 
further proved that the quotient of $AH(\pi_1(F),{\rm PSL}(2,\mathbb C))$
by ${\rm Mod}(F)$ is not even $T_1$, i.e. contains points which are not closed.

Goldman has conjectured
that ${\rm Out}(\pi_1(F))$ acts ergodically on the complement of $CC(\pi_1(F),{\rm PSL}(2,\mathbb C))$
in $X(\pi_1(F),{\rm PSL}(2,\mathbb C))$.

\medskip\noindent
{\bf Conjecture:} (Goldman) {\em If $F$ is a closed orientable surface of genus $\ge 2$, then
${\rm Out}(\pi_1(F))$ acts ergodically on $X(\pi_1(F),{\rm PSL}(2,\mathbb C))-CC(\pi_1(F),{\rm PSL}(2,\mathbb C))$.}

\medskip

The best evidence for Goldman's conjecture is provided by a result of Lee \cite{lee-Ibundle}
(see also Souto-Storm \cite{souto-storm}).

\begin{theorem}{cusps-notindd}{\rm (Lee \cite{lee-Ibundle})} If $F$ is a closed, oriented surface of genus at least two and
$U$ is an open subset of $X(\pi_1(F),{\rm PSL}(2,\mathbb C))$ which  ${\rm Out}(\pi_1(F))$
preserves and acts properly discontinuously on, then $U\cap \partial AH(\pi_1(F),{\rm PSL}(2,\mathbb C)$
is empty.
\end{theorem}

We sketch Lee's proof which makes clever use of an idea of Minsky.
First suppose that $\rho\in AH(\pi_1(F),{\rm PSL}(2,\mathbb C))$ and that there exists a simple closed curve $\alpha$ on $F$
so that $\rho(\alpha)$ is a parabolic element ($\rho$ is often referred to as a cusp). Since $\rho$ is a manifold point of
$X(\pi_1(F),{\rm PSL}(2,\mathbb C))$ (see Kapovich \cite[Thm. 8.44]{kapovich}) and (the square of the) trace of the image of
$\alpha$ is a holomorphic, non-constant function in a neighborhood of $\rho$, there exists a sequence $\{\rho_n\}$
in $X(\pi_1(F),{\rm PSL}(2,\mathbb C))$ converging to $\rho$ such that $\rho_n(\alpha)$ is a finite order elliptic for all $n$.
It follows that each $\rho_n$ is a fixed point of some power of the Dehn twist about $\alpha$. Therefore, no cusp
in $AH(\pi_1(F),{\rm PSL}(2,\mathbb C))$ can lie in any domain of discontinuity for the action of ${\rm Mod}(F)$
on $X(\pi_1(F),{\rm PSL}(2,\mathbb C))$. Since cusps
are dense in $\partial AH(\pi_1(F),{\rm PSL}(2,\mathbb C))$ (see Canary-Hersonsky \cite{canary-hersonsky}),
no representation in $\partial AH(\pi_1(F),{\rm PSL}(2,\mathbb C))$ can lie  in a domain of discontinuity.

\medskip\noindent
{\bf Remark:} Goldman's conjecture would imply that $\partial AH(\pi_1(F),{\rm PSL}(2,\mathbb C))=\partial CC(\pi_1(F),{\rm PSL}(2,\mathbb C))$
has measure zero, since it is a ${\rm Mod}(F)$-invariant subset of $X(\pi_1(F),{\rm PSL}(2,\mathbb C))$.
It would be of independent interest to prove that this is the case.

\section{Hitchin representations}

If $F$ is a closed oriented surface of genus $g\ge 2$, 
Hitchin \cite{hitchin} exhibited a component $H_n(F)$ of
$X(\pi_1(F),{\rm PSL}(n,\mathbb R))$ such that $H_n(F)$ is homeomorphic to $\mathbb R^{(n^2-1)(2g-2)}$.
Let $\tau_n:{\rm PSL}(2,\mathbb R)\to {\rm PSL}(n,\mathbb R)$ be the irreducible representation.
(The irreducible representation arises by thinking of $\mathbb R^n$ as the space of homogeneous polynomials
of degree $n-1$ in 2 variables and considering the natural action of ${\rm PSL}(2,\mathbb R)$ on this space.)
If $\rho\in{\mathcal{T}}(F)$, then we let $H_n(F)$ be the component of $X(\pi_1(F),{\rm PSL}(n,\mathbb R))$
which contains $\tau_n\circ \rho$. The subset  $\tau_n(\mathcal{T}(F))$ of $H_n(F)$ is called the {\em Fuchsian locus}.
Hitchin called $H_n(F)$ the Teichm\"uller  component because of its resemblance to Teichm\"uller space.
It is now known as the {\em Hitchin component} and representations in $H_n(F)$ are known as {\em Hitchin
representations.}
However, Hitchin's analytic techniques generated little geometric information about Hitchin representations.

Labourie \cite{labourie-anosov} introduced the dynamical notion of an Anosov representation 
(which will be further discussed in the next section) and showed
that all Hitchin representations are Anosov. As a consequence, he proved that
Hitchin  representations are discrete,  faithful and well-displacing.
Labourie \cite{labourie-ens} was then able to show that 
${\rm Mod}(F)$ acts properly discontinuously on $H_n(F)$. (Notice that 
our proof of Fricke's Theorem again generalizes to this setting, once one observes, as Labourie does,
that on any compact subset $R$ of $H_n(F)$ there exists $(J,B)$ such that if $\rho\in R$, then
$\rho$ is $(J,B)$-well-displacing.)

\begin{theorem}{}{\em (Labourie \cite{labourie-ens})}
If $F$ is a closed, oriented surface of genus $g\ge 2$ and $H_n(F)$ is the Hitchin component of 
$X(\pi_1(F),{\rm PSL}(n,\mathbb R))$, then ${\rm Mod}(F)$ acts properly discontinuously on
$H_n(F)$.
\end{theorem}

\section{Anosov representations}

In this section, we will briefly discuss  properties and examples of Anosov representations.
Anosov representations were introduced by Labourie \cite{labourie-anosov} and further
developed by Guichard and Wienhard \cite{guichard-wienhard-domains}.
For completeness, we begin with a precise definition of an Anosov representation. The casual
reader may want to simply focus on the properties in Theorems \ref{anosov-properties} and \ref{anosov-pd}
and the examples at the end of the section.  We encourage the more determined reader
to consult the original references for more details.

\medskip

Gromov  \cite{gromov} defined a geodesic
flow $(U_0\Gamma,\{\phi_t\})$ for any word hyperbolic group  $\Gamma$ (see Champetier \cite{champetier} and Mineyev
\cite{mineyev} for details). There is a natural cover $\widetilde{U_0\Gamma}$ which is identified with
$(\partial \Gamma\times \partial \Gamma-\Delta)\times\mathbb R$ where $\partial \Gamma$ is the Gromov boundary of
$\Gamma$ and $\Delta$ is the diagonal in $\partial\Gamma\times\partial \Gamma$. The geodesic flow lifts
to a flow which is simply translation in the final coordinate and $U_0\Gamma=\widetilde{U_0\Gamma}/\Gamma$.
If $\Gamma$ is the fundamental group of a closed, negatively curved manifold $N$,
then the Gromov geodesic flow $(U_0\Gamma,\{\phi_t\})$ may identified with the geodesic flow on $T^1N$.

Let $G$ be 
a semi-simple Lie group with trivial center and no compact factors,
and let $P^\pm$ be a pair of opposite proper parabolic subgroups.
Let $L=P^+\cap P^-$ be the Levi subgroup and let $S=G/L$. One may identify $S$  with an
open subset of $G/P^+\times G/P^-$,
so its tangent bundle  $TS$ admits a well-defined splitting
$TS=E^+\oplus E^-$ where $E^+|_{(x^+,x^-)}=T_{x^+}G/P^+$ and $E^-|_{(x^+,x^-)}=T_{x^-}G/P^-$. 
Let $\widetilde{F_\rho}=\widetilde{U_0\Gamma}\times S$ be the $S$-bundle over $\widetilde{U_0\Gamma}$
and let $F_\rho=(\widetilde{U_0\Gamma}\times S)/\Gamma$ be the flat $S$-bundle associated to $\rho$
(here the action on the second factor is the natural action by $\rho(\Gamma)$). The geodesic flow on 
$\widetilde{U_0\Gamma}$
extends to a flow on $\widetilde F_\rho$ which is trivial in the second factor, and descends to a flow
on $F_\rho$.  The splitting of $TS$ induces vector bundles $\tilde V_\rho^+$ and $\tilde V_\rho^-$ 
over $\widetilde F_\rho$ such that  the fibre of $\tilde V_\rho^\pm$ over the point $(z,w,x^+,x^-)\in \tilde F_\rho$
is $E^\pm|_{(x^+,x^-)}$. 
The bundles $\tilde V_\rho^+$ and $\tilde V_\rho^-$ descend to bundles $V_\rho^+$ and $V_\rho^-$ over $F_\rho$
and there is again an induced flow covering the geodesic flow.

A representation $\rho:\Gamma\to G$ is {\em $P^\pm$-Anosov } if there exists a  section 
\hbox{$\sigma: U_0\Gamma\to F_\rho$} which is flat along flow lines
so that the induced flow on the bundle $\sigma^*V^-$ is contracting
(i.e. given any metric $||\cdot||$ on $\sigma^*V^-$ there exists $t_0>0$ so that if $v\in \sigma^*V^-$,
then $||\phi_{t_0}(v)||\le \frac{1}{2} ||v||$) and 
the induced flow on the bundle $\sigma^*V^+$ is dilating (i.e. the inverse flow is
contracting). 
In particular, $\sigma(U_0\Gamma)$ is a hyperbolic set for the lift of the geodesic flow to $F_\rho$.
We recall that the section is said be flat along flow lines
if it lifts to a section $\tilde\sigma$ of $\widetilde F_\rho$ whose second factor is constant along flow lines.

\medskip

Let ${\rm Anosov}_{P^\pm}(\Gamma,G)\subset X(\Gamma,G)$ be the set of 
(conjugacy classes of) \hbox{$P^\pm$-Anosov} representations of $\Gamma$ into $G$.
We further let  ${\rm Anosov}(\Gamma,G)\subset X(\Gamma,G)$
denote the set of  (conjugacy classes of ) representations of $\Gamma$ into $G$
which are Anosov with respect to some pair of opposite proper parabolic subgroups of $G$.

The following properties are established by Labourie \cite{labourie-anosov,labourie-ens} and 
Guichard-Wienhard \cite{guichard-wienhard-domains}.

\begin{theorem}{}{\em (Labourie \cite{labourie-anosov,labourie-ens},
Guichard-Wienhard \cite{guichard-wienhard-domains})}
\label{anosov-properties}
If $\Gamma$ is a torsion-free word hyperbolic group, $G$ is
a semi-simple Lie group with trivial center and no compact factors,
and $P^\pm$ is a pair of opposite proper parabolic subgroups, then
\begin{enumerate}
\item
${\rm Anosov}_{P^\pm}(\Gamma,G)$ is an open subset of $X(\Gamma,G)$.
\item
If $\rho\in{\rm Anosov}_{P^\pm}(\Gamma,G)$, then $\rho$ is discrete and faithful.
\item
If $\rho\in{\rm Anosov}_{P^\pm}(\Gamma,G)$, then $\rho$ is well-displacing.
\item
If $\rho\in{\rm Anosov}_{P^\pm}(\Gamma,G)$ and $x\in X$, then  the orbit map $\tau_{\rho,x}:\Gamma\to X$ is a quasi-isometric embedding.
\end{enumerate}
\end{theorem}

Once one has established this result, one may use the proof of Fricke's Theorem 
given in section 2 to obtain:

\begin{theorem}{anosov-pd}{\em (Guichard-Wienhard \cite[Cor. 5.4]{guichard-wienhard-domains})}
If $\Gamma$ is a torsion-free word hyperbolic group and $G$ is
a semi-simple Lie group with trivial center and no compact factors,
then ${\rm Out}(\Gamma)$ acts properly discontinuously on
${\rm Anosov}(\Gamma,G)$.
\end{theorem}

\noindent
{\bf Remark:} Guichard and Wienhard only state Theorem \ref{anosov-pd} for surface groups
and free groups, but the only added ingredient needed for the full statement we give is
Proposition \ref{gens-suffice}.

\medskip\noindent
{\bf Examples of Anosov representations:}
(1) If $B^+$ is the subgroup of upper triangular matrices in ${\rm PSL}(n,\mathbb R)$ and
$B^-$ is the subgroup of lower triangular matrices, then $B^\pm$ is a pair of opposite
parabolic subgroups of ${\rm PSL}(n,\mathbb R)$. Labourie \cite{labourie-anosov}
showed that Hitchin representations are \hbox{$B^\pm$-Anosov.}

\medskip

(2) Suppose that  $G={\rm PSL}(n,\mathbb R)$, $P^+$ is the stabilizer of a line and $P^-$ is
the stabilizer of an orthogonal hyperplane. Then $P^\pm$ is a pair of opposite parabolic subgroups.
Representations into ${\rm PSL}(n,\mathbb R)$ which are $P^\pm$-Anosov are called
{\em convex Anosov}.

If $\rho:\Gamma\to {\rm PSL}(n,\mathbb R)$, $\Omega$ is a strictly convex domain in $\mathbb P(\mathbb R^n)$ and 
$\rho(\Gamma)$ acts properly discontinuously and cocompactly on $\Omega$, then it follows
immediately from work of Benoist \cite{benoist-divisible1}, that $\rho$ is convex Anosov. Such
representations are called {\em Benoist representations}. Benoist \cite{benoist-divisible3}
further proved that if a component $C$ of $X(\Gamma,{\rm PSL}(n,\mathbb R))$ contains
a Benoist representation, then it consists entirely of Benoist representations. We call such components
{\em Benoist components.}

(3) If $G$ is a rank one semisimple Lie group, then $\rho:\Gamma\to G$ is Anosov if and only
$\rho$ is convex cocompact (see Guichard-Wienhard \cite[Thm. 5.15]{guichard-wienhard-domains}).
For this reason, along with the properties in Theorem 
\ref{anosov-properties},  it is natural to think of Anosov representations as the higher rank analogues
of convex cocompact representations into rank one Lie groups. 

We recall that Kleiner-Leeb \cite{kleiner-leeb}
and Quint \cite{quint} proved that if $G$ has  rank at least 2, then every convex compact subgroup of $G$
is, up to finite index, a product of cocompact lattices and convex cocompact subgroups of rank one Lie groups.
Therefore, convex cocompact representations will not yield a robust class of representations in higher rank.

\section{${\rm PSL}(2,\mathbb C)$-character varieties}

In this section, we will restrict to the setting  where $G={\rm PSL}(2,\mathbb C)$, $M$ is a compact, orientable
3-manifold with non-empty boundary, such that no component of $\partial M$ is a torus, and $\Gamma=\pi_1(M)$.
In light of Goldman's conjecture it is natural to ask whether or not
$CC(\pi_1(M),{\rm PSL}(2,\mathbb C))$ is a maximal domain of discontinuity for the action of
${\rm Out}(\pi_1(M))$ on $X(\pi_1(M),{\rm PSL}(2,\mathbb C))$.
We will see that very often 
$CC(\pi_1(M),{\rm PSL}(2,\mathbb C))$ is not a maximal domain of discontinuity. It appears
that the dynamics in
the case where $M=F\times [0,1]$ are quite different than the dynamics in the general case.

\subsection{Free groups}
\label{freepsl2c}

Minsky \cite{minsky-PS} first studied the case where $M$ is a handlebody of genus $n\ge 2$
and $\Gamma$ is the free group $F_n$ of rank $n$. He exhibits a domain of discontinuity $PS(F_n)$ for
the action of ${\rm Out}(F_n)$ on $X(F_n,{\rm PSL}(2,\mathbb C))$ which is strictly larger than 
$CC(F_n,{\rm PSL}(2,\mathbb C))$. We note that $CC(F_n,{\rm PSL}(2,\mathbb C))$ is the space of Schottky
groups of genus $n$ and is an open submanifold of $X(F_n,{\rm PSL}(2,\mathbb C))$ with infinitely 
generated fundamental group (see Canary-McCullough \cite[Chap. 7]{canary-mccullough}
for a general discussion of the structure
of $CC(\Gamma, {\rm PSL}(2,\mathbb C))$).

We recall that the Cayley
graph  $C_{F_n}$ of $F_n$ with respect to a minimal set of generators is a  $2n$-valent tree. 
We say that  \hbox{$\rho\in X(F_n,{\rm PSL}(2,\mathbb C))$}
is  {\em primitive-stable} if and only if for some $x\in\mathbb H^3$ there exists $(K,C)$ such that the restriction of
the extended orbit map $\bar\tau_{\rho,x}:C_{F_n}\to \mathbb H^3$ 
to the axis of any primitive element in $C_{F_n}$ is a $(K,C)$-quasi-isometric embedding.
We recall that an element $g\in F_n$ is {\em primitive} if it is a member of a minimal generating set for $F_n$.

Let $PS(F_n)\subset X(F_n,{\rm PSL}(2,\mathbb C))$ be the set of (conjugacy classes of) primitive-stable
representations.

\begin{theorem}{minsky-psthm}{\rm (Minsky \cite{minsky-PS})}
If $n\ge 2$, then $PS(F_n)$ is an open subset of $X(F_n,{\rm PSL}(2,\mathbb C))$ containing
$CC(\Gamma,{\rm PSL}(2,\mathbb C))$ such that
\begin{enumerate}
\item
${\rm Out}(F_n)$ preserves and acts properly discontinuously on $PS(F_n)$,
\item
$PS(F_n)$ intersects $\partial CC(F_n,{\rm PSL}(2,\mathbb C))$, and
\item
$PS(F_n)$ contains representations which are indiscrete and not faithful.
\end{enumerate}
\end{theorem}

The proof that $PS(F_n)$ is open is much the same as the proof of Proposition \ref{openinrank1}.

The proof that ${\rm Out}(F_n)$ acts properly discontinuously on $PS(F_n)$ again follows the outline
of our proof of  Fricke's Theorem. 
One first proves, just as in the proof of Proposition \ref{qigiveswd}, that  if $\rho\in PS(F_n)$, then $\rho$ is well-displacing on primitive elements,
i.e. that there exists $J$ and $B$ such that if $\gamma\in F_n$ is primitive, then
$$\frac{1}{J}||\gamma||-B\le \ell_\rho(\gamma)\le J||\gamma||+B.$$
Notice  that if $\alpha$ and $\gamma$ are distinct elements
of a minimal generating set for $F_n$, then $\alpha\gamma$ is primitive, so, by Proposition \ref{gens-suffice}, if
$\{\phi_n\}$ is a sequence of distinct elements of ${\rm Out}(F_n)$,  then, up to subsequence, there exists a primitive element
$\beta\in \Gamma$, so that $||\phi_n^{-1}(\beta)||\to\infty$. The remainder of the proof is the same.

The most significant, and difficult, part of the proof of Theorem \ref{minsky-psthm} is
the proof of item (2). Minsky makes clever use of Whitehead's algorithm  to determine whether or
not an element of a free group is primitive, to prove that
if $\rho$ is geometrically finite and $\rho(\gamma)$ is parabolic for some curve $\gamma$ in the Masur
domain, then $\rho$ is primitive-stable. Item (3) then follows nearly immediately from (2).
A more thorough discussion of this portion of the proof
is outside the parameters of this survey article, but we encourage reader to consult Minsky's
beautiful paper.

\medskip

If $n\ge 3$,  then Gelander  and Minsky defined 
an  open subset $R(F_n)$  of $X(F_n,{\rm PSL}(2,\mathbb C))$ which ${\rm Out}(F_n)$
acts ergodically on. A representation \hbox{$\rho:F_n\to {\rm PSL}(2,\mathbb C)$} is
{\em redundant} if and only if there exists a non-trivial free decomposition
$F_n=A*B$  so that $\rho(A)$ is dense in ${\rm PSL}(2,\mathbb C)$.

\begin{theorem}{}{\em (Gelander-Minsky \cite{gelander-minsky})}
If $n\ge 3$, then  the set $R(F_n)$ of redundant representations
of $F_n$ into ${\rm PSL}(2,\mathbb C)$ is an open subset of $X(F_n,{\rm PSL}(2,\mathbb C))$ and
${\rm Out}(F_n)$ acts ergodically on $R(F_n)$.
\end{theorem}

It is easy to see that no redundant representation can lie in a domain of discontinuity for
the action of ${\rm Out}(F_n)$. Suppose that $\rho$ is redundant, that $F_n=A*B$, $\rho(A)$
is dense in ${\rm PSL}(2,\mathbb C)$ and $B=<b>$. Let $\{ a_n\}$ be a sequence of 
(non-trivial) elements
of $A$ so that $\{\rho(a_n)\}$ converges to the identity. For each $n$, let $\phi_n\in {\rm Out}(F_n)$
be the automorphism which is trivial on $A$ and takes $b$ to $a_nb$. Then $\{\rho_n\circ\phi_n\}$ converges
to $\rho$.

\medskip

Gelander and Minsky \cite{gelander-minsky} asked the following natural question:

\medskip\noindent
{\bf Question:} {\em If $n\ge 3$, does $PS(F_n)\cup R(F_n)$ have full measure in $X(F_n,{\rm PSL}(2,\mathbb C))$?}

\medskip

A positive answer to this question would give a very satisfying dynamical dichotomy in
the case of the  ${\rm PSL}(2,\mathbb C)$-character variety of a free group.

\medskip

Lee \cite{lee-compbody} extend Minsky's results into the setting of fundamental
groups of compression bodies (without toroidal boundary components), in which case
$\pi_1(M)$ is the free product of a free group and a finite number of closed orientable surface groups.

\begin{theorem}{comp-ss}{\rm (Lee \cite{minsky-PS})}
If $\Gamma=\Gamma_1*\cdots*\Gamma_n$, $n\ge 2$ and each $\Gamma_i$ is either infinite
cyclic or isomorphic to the fundamental group of a closed orientable surface of genus at least 2,
then  there exists an open subset $SS(\Gamma)$ of $X(\Gamma,{\rm PSL}(2,\mathbb C))$
containing $CC(\Gamma,{\rm PSL}(2,\mathbb C))$ such that
\begin{enumerate}
\item
${\rm Out}(\Gamma)$ preserves and acts properly discontinuously on $SS(\Gamma)$,
\item
$SS(\Gamma)$ intersects $\partial CC(\Gamma,{\rm PSL}(2,\mathbb C))$, and
\item
$SS(\Gamma)$ contains representations which are indiscrete and not faithful.
\end{enumerate}
\end{theorem}

If $\Gamma$ is as in Theorem \ref{comp-ss}, but is not a free product of two surface groups,
then Lee defines a representation $\rho\in X(\Gamma,{\rm PSL}(2,\mathbb C))$ 
to be {\em separable-stable} if there exists $(K,C)$ and $x\in\mathbb H^3$ so that the restriction 
of the extended orbit map $\tau_{\rho,x}:C_\Gamma\to \mathbb H^3$ to any geodesic joining
the fixed points of a separable element of $\Gamma$ is a $(K,C)$-quasi-isometric embedding.
(Recall that $\gamma\in\Gamma$ is said to be {\em separable} if it lies in a factor of a non-trivial
free decomposition of $\Gamma$.) With this definition the open-ness of the set $SS(\Gamma)$ 
of separable-stable representations and the proper discontinuity of the action of ${\rm Out}(\Gamma)$
follows much as in the proofs of Fricke's Theorem and Minsky's Theorem \ref{minsky-psthm}.

Again, the proof of (2) is the most significant, and difficult, portion of the proof of Theorem \ref{comp-ss}.
Lee makes use of a generalization of Whitehead's algorithm, due
to Otal \cite{otal-thesis}, and there are significant new technical issues in generalizing Minsky's
arguments to this setting.

\medskip\noindent
{\bf Remarks:}
(1) Tan, Wong and Zhang \cite{TWZ}, in work which predates Minsky's work, 
defined an open subset $BQ$ of $X(F_2,{\rm SL}(2,\mathbb C))$
which contains all lifts of convex cocompact representations of $F_2$ into ${\rm PSL}(2,\mathbb C)$
such that ${\rm Out}(F_2)={\rm GL}(2,\mathbb Z)$ acts properly discontinuously on $BQ$.
(The representations in
$BQ$ are those that satisfy the Bowditch $Q$-conditions, which Bowditch originally defined in 
\cite{bowditch}.) It is natural to ask whether
$BQ$ is a maximal domain of discontinuity for the action of ${\rm Out}(F_2)$.

The set $BQ$ also contains all lifts of primitive-stable representations. It is an open
question  whether or not $BQ$ consists entirely of lifts of primitive-stable representations.

(2) Jeon, Kim, Lecuire and Ohshika \cite{JKLO} give a complete characterization of which points
in $AH(F_n, {\rm PSL}(2,\mathbb C))$ are primitive-stable.
 
(3) Minsky and Moriah \cite{minsky-moriah} have exhibited large classes of non-faithful primitive-stable
representations whose images are lattices in ${\rm PSL}(2,\mathbb C)$.

\subsection{Freely indecomposable groups}

In this section, we consider the case that $M$ is a compact, orientable \hbox{3-manifold}
whose interior admits a complete hyperbolic metric  such that $\partial M$ contains
no tori and $\pi_1(M)$ is freely
indecomposable, i.e. does not admit a free decomposition. Equivalently,
we could require that $M$ is irreducible (i.e. every embedded sphere in $M$ bounds a ball in $M$)
has incompressible boundary (i.e. if $S$ is a component of $\partial M$, then
$\pi_1(S)$ injects into $\pi_1(M)$) and $\pi_1(M)$ is infinite and contains no rank two free abelian subgroups.

One may combine work of Canary-Storm \cite{canary-storm2}  and Lee \cite{lee-Ibundle} to obtain

\begin{theorem}{biggerDD}{\em (Canary-Storm \cite{canary-storm2}, Lee \cite{lee-Ibundle})}
If $M$ is a compact, orientable 3-manifold
whose interior admits a complete hyperbolic metric  such that $\partial M$ contains
no tori and $\pi_1(M)$ is freely
indecomposable, and $M$ is not homeomorphic to 
a trivial interval bundle $F\times [0,1]$,
then there exists an open subset $W(M)$ of $X(\pi_1(M),{\rm PSL}(2,\mathbb C))$ such that
\begin{enumerate}
\item
$W(M)$ is invariant under ${\rm Out}(\pi_1(M))$, 
\item ${\rm Out}(\pi_1(M))$ acts properly discontinuously on $W(M)$,
\item $CC(\pi_1(M),{\rm PSL}(2,\mathbb C))$ is a proper subset of $W(M)$.
\item
If $\rho\in AH(\pi_1(M)),{\rm PSL}(2,\mathbb C))$  and $\rho(\pi_1(M))$ contains no
parabolic elements, then  $\rho\in W(M)$. 
\item $W(M)$ contains representations with indiscrete image which are not faithful.
\end{enumerate}
\end{theorem}

Lee \cite{lee-Ibundle} handles the case where $M$ is a twisted interval bundle (i.e. an interval bundle which is
not a product). In this case, $\pi_1(M)$ is the fundamental group of a non-orientable surface, so
we see that there is even a substantial difference between the dynamics on 
the ${\rm PSL}(2,\mathbb C)$-character varieties associated
to orientable and non-orientable surface groups. Lee \cite{lee-Ibundle} also characterizes exactly which points in
$AH(\pi_1(M),{\rm PSL}(2,\mathbb C))$
can lie in a domain of discontinuity for the action of ${\rm Out}(\pi_1(M))$ when $M$ is a 
twisted interval bundle.

Canary and Storm \cite{canary-storm2} also  determine exactly when $AH(\pi_1(M),{\rm PSL}(2,\mathbb C))$
is entirely contained in a domain of discontinuity.  We recall that
an {\em essential annulus} in $M$ is an embedded annulus
$A\subset M$ such that $\partial A\subset \partial M$, $\pi_1(A)$
injects into $\pi_1(M)$ and $A$ cannot be homotoped (rel $\partial A$)
into $\partial M$.
We say that
an essential annulus is {\em primitive} if $\pi_1(A)$ is a maximal abelian
subgroup of $\pi_1(M)$.
Johannson \cite{johannson} proved that  ${\rm Out}(\pi_1(M))$ is finite if and only if $M$ contains no essential annuli.
So, ${\rm Out}(\pi_1(M))$ acts properly discontinuously on $X(\pi_1(M),{\rm PSL}(2,\mathbb C))$ if and only if 
$M$ contains no essential annuli.

Canary and Storm \cite{canary-storm2} show that if $M$ contains a primitive essential annulus $A$,
then there exists $\rho\in AH(\pi_1(M),{\rm PSL}(2,\mathbb C))$ which cannot be in any domain of discontinuity for the 
action of ${\rm Out}(\pi_1(M))$ on $X(\pi_1(M), {\rm PSL}(2,\mathbb C))$. One can use a generalization of
Lee's proof of Theorem \ref{cusps-notindd} to establish this.
Thurston's geometrization theorem for pared manifolds (see Morgan \cite{morgan})
implies that there exists $\rho\in AH(\pi_1(M),{\rm PSL}(2,\mathbb C))$ 
such that $\rho(\pi_1(A))$ consists of parabolic elements.
Then one can again show that
$\rho$ is the limit of a sequence $\{\rho_n\}$ of representations in $X(\pi_1(M),{\rm PSL}(2,\mathbb C))$ such that,
for all $n$, $\rho_n(\pi_1(A))$ is a finite order elliptic group, so that $\rho_n$ is fixed by some power of
the Dehn twist about $A$.  They also give a more complicated argument  which shows
that there exists $\rho\in AH(\pi_1(M),{\rm PSL}(2,\mathbb C))$ which cannot be in any domain of discontinuity for the 
action of ${\rm Out}(\pi_1(M))$ on $AH(\pi_1(M), {\rm PSL}(2,\mathbb C))$.

The presence of primitive essential annuli
is the only obstruction to $AH(\pi_1(M),{\rm PSL}(2,\mathbb C))$ lying entirely in a domain of discontinuity.

\begin{theorem}{AHinDD}{\em (Canary-Storm \cite{canary-storm2})}
If $M$ is a compact, orientable 3-manifold
whose interior admits a complete hyperbolic metric  such that $\partial M$ contains
no tori and $\pi_1(M)$ is freely
indecomposable, then
${\rm Out}(\pi_1(M))$ act properly discontinuously on an open neighborhood of 
$AH(\pi_1(M),{\rm PSL}(2,\mathbb C))$ in $X(\pi_1(M),{\rm PSL}(2,\mathbb C))$ if and only if $M$ 
contains no primitive essential annuli. 
\end{theorem}

The techniques used by Canary and Storm \cite{canary-storm2}  to handle the cases
in Theorems \ref{biggerDD} and \ref{AHinDD} when $M$ is not an interval bundle are substantially
different than the technique used by Minsky and Lee.
Canary and Storm make crucial use of the study of mapping class groups of
3-manifolds developed by Johannson \cite{johannson} and extended by McCullough \cite{mccullough-VGF}
and Canary-McCullough \cite{canary-mccullough}. 

We provide a brief sketch of the techniques involved.
The rough idea is that there is a finite index subgroup ${\rm Out}_0(\pi_1(M))$ of ${\rm Out}(\pi_1(M))$ which is constructed
from subgroups generated by Dehn twists about essential annuli in $M$ and mapping class groups of base surfaces
of interval bundle components of the characteristic submanifold. It turns out that it  suffices
to require that there are subgroups ``registering'' each of these two types of building blocks for ${\rm Out}(\pi_1(M))$
such that  the restriction of the representation to each such subgroup is 
convex cocompact  (or even just primitive-stable). A more detailed sketch follows.

Suppose that $M$ is a compact 3-manifold
whose interior admits a complete hyperbolic metric  such that $\partial M$ contains
no tori, $\pi_1(M)$ is freely
indecomposable and that $M$ is not an interval bundle.
Then $M$ contains a submanifold $\Sigma$, called the
{\em characteristic submanifold},
such that each component of $\Sigma$ is either an $I$-bundle which intersects $\partial M$
in its associated $\partial I$-bundle or a solid torus. Moreover, each component of the frontier
of $\Sigma$ is an essential annulus in $M$ and every essential annulus in $M$ is properly isotopic
into $\Sigma$. (The characteristic submanifold was developed by Jaco-Shalen 
\cite{jaco-shalen} and Johannson \cite{johannson}. For a discussion of the characteristic manifold
in the hyperbolic setting, see Morgan \cite{morgan} or Canary-McCullough \cite{canary-mccullough}.)

Canary and McCullough \cite{canary-mccullough} showed that, in this setting, there exists a finite
index subgroup ${\rm Out}_1(\pi_1(M))$ of ${\rm Out}(\pi_1(M))$ consisting of outer automorphisms
which are induced by homeomorphisms of $M$. McCullough \cite{mccullough-VGF} showed
that there is a further finite index subgroup ${\rm Out}_0(\pi_1(M))$  consisting of outer automorphisms induced by
homeomorphisms fixing each component of $M-\Sigma$ and a short exact sequence
$$1\longrightarrow K(M)\longrightarrow {\rm Out}_0(\pi_1(M))\longrightarrow \oplus {\rm Mod}_0(F_i)\longrightarrow 1$$
where 
$K(M)$ is the free abelian group generated by Dehn twists about essential annuli in the frontier of $\Sigma$,
$\{ F_1,\ldots, F_k\}$ are the base surfaces of the interval bundle components
and ${\rm Mod_0}(F_i)$ denotes the group of isotopy classes of homeomorphisms of $F_i$ whose
restriction to each boundary component is isotopic to the identity.
(Sela \cite{sela} used his JSJ splitting of a torsion-free word hyperbolic group 
to obtain a version of this short exact sequence in the
setting of outer automorphism groups of torsion-free word hyperbolic groups, see also Levitt \cite{levitt}.)

A {\em characteristic collection of annuli} for $M$ is either (a) the collection of components of the
frontier of a solid torus component of $\Sigma$, or (b) a component of the frontier of an interval
bundle component of $\Sigma$ which is not isotopic into any solid torus component of $\Sigma$.
Let $\{{\mathcal{A}}_1,\ldots,{\mathcal{A}}_r\}$ be the set of all characteristic collections of annuli.
If $K_i$ is the free abelian group generated by Dehn twists in the elements of ${\mathcal{A}}_i$, then
$K(M)=\oplus K_i$. 

If ${\mathcal{A}}_i$ is  a characteristic collection of annuli of type (a), 
we say that a free subgroup $H$ of $\pi_1(M)$ {\em registers} ${\mathcal{A}}_i=\{ A_{i1},\ldots,A_{is}\}$ if,
it is freely generated by loops $\{ \alpha_1,\alpha_0,\ldots,\alpha_s\}$ so that $\alpha_0$ is the
core curve of the solid torus component of $\Sigma$ and each $\alpha_j$,
for $j=1,\ldots,s$, is a loop based at $x_0\in\alpha_0$ which intersection $A_{ij}$ exactly twice and intersects no other
annulus in any characteristic collection of annuli for $M$. (We give a similar definition for characteristic
collections of annuli of type (b).) Then $K_i$ preserves the subgroup $H_i$ and injects into ${\rm Out}(H_i)$.

We then say that  a representation $\rho\in X(\pi_1(M),{\rm PSL}(2,\mathbb C))$ lies in $W(M)$
if 
\begin{enumerate}
\item for every characteristic collection of annuli ${\mathcal{A}}_i$ there
is a registering subgroup $H_i$ such that $\rho|_{H_i}$ is primitive-stable, and 
\item for every 
interval bundle component $\Sigma_j$ of $\Sigma$, whose base surface is not a two-holed
projective plane or three-holed sphere,
$\rho|_{\pi_1(\Sigma_j)}$ is primitive stable. 
\end{enumerate}

Suppose that $\rho\in W(M)$.
If $\{\phi_n\}$ is a sequence in $K(M)$, which projects
to a sequence $\{\psi_n\}$ of distinct elements of $K_i$. Then, if $H_i$ is a registering subgroup for $\mathcal{A}_i$
such that $\rho|_{H_i}$ is primitive-stable, then, by Minsky's Theorem,
\hbox{$\rho|_{H_i}\circ \psi_n^{-1}\to\infty$}
in $X(H_i,{\rm PSL}(2,\mathbb C))$. Therefore, \hbox{$\rho\circ \phi_n^{-1}\to\infty$} in
$X(\pi_1(M)),{\rm PSL}(2,\mathbb C))$.
If $\{\phi_n\}$ is a sequence in ${\rm Out}_1(\pi_1(M))$ whose projection to ${\rm Mod}_0(F_j)$ 
produces a sequence $\{\psi_n\}$ of distinct elements, then,
since $\rho|_{\pi_1(\Sigma_j)}$ is primitive-stable,
\hbox{$\rho|_{\pi_1(\Sigma_j)}\circ \psi_n^{-1}\to\infty$}
in $X(\pi_1(\Sigma_j),{\rm PSL}(2,\mathbb C))$. Thus, \hbox{$\rho\circ \phi_n^{-1}\to\infty$} in
$X(\pi_1(M)),{\rm PSL}(2,\mathbb C))$. 
Combining these two cases, and doing a little book-keeping
yields a proof of item (2) in Theorem \ref{biggerDD}. If $\rho\in AH(\pi_1(M),{\rm PSL}(2,\mathbb C))$
and $\rho(\pi_1(\Sigma))$ consists entirely of hyperbolic elements, then an application of
the ping-pong lemma and the covering theorem guarantee that $\rho\in W(M)$, which
may be used to establish  items (3) and (4) in Theorem \ref{biggerDD}. 

If $M$ contains no primitive essential annuli and $\rho\in AH(\pi_1(M),{\rm PSL}(2,\mathbb C))$,
then $\rho(\pi_1(\Sigma))$ cannot
contain any parabolic elements, which do not lie in the image of an interval bundle component of
$\Sigma$ whose base surface is a two-holed projective plane. So, one can again use
the ping-pong lemma and the covering theorem to show that  $\rho\in AH(M)$. Therefore, if $M$ contains
no primitive essential annuli, then
$AH(\pi_1(M),{\rm PSL}(2,\mathbb C))$ is contained in $W(M)$, which is an open set
which is a domain of discontinuity for the action of ${\rm Out}(\pi_1(M))$ on $X(\pi_1(M),{\rm PSL}(2,\mathbb C))$.
Since we have previously noted that if $M$ contains a primitive essential annulus, then no
domain of discontinuity can contain $AH(\pi_1(M),{\rm PSL}(2,\mathbb C))$, this completes
the proof of Theorem \ref{AHinDD}.

\medskip\noindent
{\bf Remarks:} (1)
Theorems \ref{biggerDD} and \ref{AHinDD} were generalized to the case where $\partial M$ is allowed
to contain tori by Canary and Magid \cite{canary-magid}.

(2) If $M$ is a compact, orientable 3-manifold whose interior
admits a hyperbolic metric, then  Canary and Storm \cite{canary-storm2} study the quotient  moduli space
$$AI(M)=AH(\pi_1(M),{\rm PSL}(2,C))/{\rm Out}(\pi_1(M)).$$
They show that $AI(M)$ 
is $T_1$ (i.e. all points are closed) unless $M$ is a product interval bundle.
However, if $M$ contains a primitive essential annulus, then $AI(M)$
is not Hausdorff. If $M$ contains no primitive annuli and has no toroidal boundary components, they
show that $AI(M)$ is Hausdorff.

\section{Appendix}

In this appendix, we give a proof of Proposition \ref{gens-suffice}.

\medskip\noindent{\bf Proposition \ref{gens-suffice}.} {\em
If $\Gamma$ is a torsion-free word hyperbolic group with finite generating set $S$,
then there exists a finite set $\mathcal B$ of elements
of $\Gamma$ such that for any $K$,
$$\{\phi\in {\rm Out}(\Gamma)\ |\  ||\phi(\beta)||\le K \ \ {\rm if}\ \beta\in{\mathcal{B}} \}$$
is finite.

Moreover, if $\Gamma$ admits a convex cocompact action on the symmetric space associated
to a rank one Lie group, then one may take
$\mathcal B$ to consist of all the elements of $S$ and all products of two elements of $S$.
}

\medskip\noindent
{\em Proof:}
We first prove our result in the case that  $\Gamma$ acts convex cocompactly on the symmetric space $X$
associated to a rank one Lie group $G$, i.e. that there exists a convex cocompact
representation $\rho:\Gamma\to G$. We will give a proof in this case which readily generalizes to
the general hyperbolic setting.

We first show that if $||\phi(\beta)||\le K$ for all $\beta\in\mathcal{B}$, then there is an upper
bound $C_1$ on the distance between the axes of the images of the generators. We next show
that there exists a point which lies within $M$ of the axis of the image of every generator. Since,
by conjugating the representative of the automorphism,
we may assume that this point lies in a compact set, the image of each generator is
determined up to finite ambiguity, which implies that the outer automorphism is determined up
to finite ambiguity.

\medskip

Since $\rho$ is convex cocompact, 
there exists $(J,B)$ so that $\rho$ is $(J,B)$-well-displacing.
Moreover, there exists
$\eta>0$ so that $\ell_\rho(\gamma)>\eta$ if $\gamma\in \Gamma-\{id\}$.

Suppose that $S=\{\alpha_1,\ldots,\alpha_n\}$ is a generating set for $\Gamma$,
$\phi\in {\rm Out}(\Gamma)$ and $||\phi(\gamma)||\le K$ if $\gamma$ is an element of $S$ or
a product of two elements of $S$. Let $\bar\phi$ be a representative of the conjugacy class $\phi$.

For each $i$, let $A_i$ be the axis of $\rho(\bar\phi(\alpha_i))$ and 
let $p_i:X\to A_i$ be the nearest
point projection onto $A_i$. We first prove that there exists $C_1$, depending only
on $X$, $K$, $J$, $B$ and $\eta$, so that
$$d(A_i,A_j)\le C_1$$
for all $i,j\in\{1,\ldots,n\}$.

For each distinct pair $i,j\in\{1,\ldots,n\}$, choose $x_{ij}\in A_i$ and $y_{ij}\in A_j$ so that
$\overline{x_{ij}y_{ij}}$ is the unique common perpendicular joining $A_i$ to $A_j$.
(If $A_i$ and $A_j$ intersect, let $x_{ij}=y_{ij}$ be the point of intersection).
Let $p_{ij}:X\to \overline{x_{ij}y_{ij}}$ be the nearest point projection onto $\overline{x_{ij}y_{ij}}$.
Notice that $p_{ij}(A_i)=x_{ij}$ and $p_i(\overline{x_{ij}y_{ij}})=x_{ij}$.
Let $x_{ij}^\pm$ be the points on $A_i$ which are a distance $\eta\over 2$ from
$x_{ij}$. Let $Q_{ij}=p_i^{-1}(\overline{x_i^-x_i^+})$ where $\overline{x_i^-x_i^+}$ is the geodesic segment joining
$x_i^-$ to $x_i^+$. There exists $C_0>0$ so that
$$p_{ij}(X-Q_{ij})\subset B(x_{ij},C_0).$$
(Since $X$ is $CAT(-1)$, the constant $C_0$ which works for $\mathbb H^2$ will also work for $X$.)

Let $R_{ij}=p_j^{-1}(\overline{x_{ji}^+x_{ji}^-})=Q_{ji}$. Then
$$p_{ij}(X-R_{ij})\subset B(y_{ij},C_0).$$
So, if 
$$d(A_i,A_j)>C_1=2C_0+JK+B,$$
then, since $p_{ij}$ is 1-Lipschitz,
$X-Q_{ij}$ and $X-R_{ij}$ are disjoint and 
$$d(X-Q_{ij},X-R_{ij})>d(A_i,A_j)-2C_0>JK+B.$$
Since $Q_{ij}$ is contained in a fundamental domain for the action of $<\rho(\bar\phi(\alpha_i))>$ and $R_{ij}$ is
contained in a fundamental  domain for the action of $<\rho(\bar\phi(\alpha_j))>$,
$Q_{ij}\cap R_{ij}$ is contained in a fundamental domain for the action of the  subgroup of $\rho(\Gamma)$
generated by $\rho(\bar\phi(\alpha_i))$ and $\rho(\bar\phi(\alpha_j))$. Moreover, $Q_{ij}\cap R_{ij}$ is
contained in a fundamental domain for $\rho(\bar\phi(\alpha_i\alpha_j))$ and the fixed points of $\rho(\bar\phi(\alpha_i\alpha_j))$
are contained in $\overline{X-Q_{ij}}$ and $\overline{X-R_{ij}}$.
Since $d(X-Q_{ij},X-R_{ij})>JK+B$,
$\ell_\rho(\bar\phi(\alpha_i\alpha_j))>JK+B$.  Since $\rho$ is $(J,B)$-well-displacing,
$||\bar\phi(\alpha_i\alpha_j)||>K$. Therefore, we may conclude that
$$d(A_i,A_j)\le C_1.$$

We next claim that  there exists $M>0$ and a
point $x\in X$ so that $d(x,A_i)\le M$ for all $i=1,\ldots, n$. 
We first prove the claim for each collection of 3 generators.
For any distinct $i,j,k$ in $\{1,\ldots,n\}$, we may
construct a geodesic hexagon consisting of
the unique common perpendiculars joining $A_i$ to $A_j$, $A_j$ to $A_k$ and $A_k$ to $A_i$ and
the geodesic segments of $A_i$, $A_j$ and $A_k$ joining them. (If $A_l$ and $A_m$ intersect, 
replace the unique common perpendicular with the point of intersection and regard it as a degenerate
edge of the hexagon.) Since $X$ is $CAT(-1)$, there exists $W>0$ such that every geodesic triangle
contains a ``central point'' which lies within $W$ of each side of the triangle.
Let  $z_{ijk}$ be a``central'' point of a geodesic triangle  $T_{ijk}$ joining alternate endpoints of the hexagon,
so $z_{ijk}$ lies within $W$ of each edge of $T_{ijk}$.
If $s\in\{i,j,k\}$, then there exists a triangle $T_s$ which shares an edge with $T_{ijk}$ and whose
other sides  are  the edge  of the hexagon contained in $A_s$ and one of the common perpendiculars.
Since this triangle is $\delta$-slim (for some $\delta$ depending only on $X$)
and each common perpendicular has length at most $C_1$, we see that
$$d(z_{ijk},A_s)\le M_3= W+C_1+\delta$$
for all $s\in\{i,j,k\}$.

Now notice that there is an upper bound $D$ on the
diameter of $p_i(A_j)$, for all $i$ and $j$, which depends only on $K$ and 
$\rho(\Gamma)$ (since, up to conjugacy,
there are only finitely  many pairs of  elements of $\rho(\Gamma)$ with translation distance at most $JK+B$ whose
axes are separated by at most $C_1$.) Therefore,
$$d(p_1(z_{123}),p_1(z_{i2k}))\le 2M_3+D $$
if $k\in\{4,\ldots,n\}$
(since both lie within $M_3$ of $p_1(A_2)$). It follows, since $d(z_{12k},A_k)<M_3$, that if $x=z_{123}$,
then
$$d(x,A_k)\le M=4M_3+D$$
for all $k\in\{1,\ldots,n\}$.

We may then assume  that the
representative $\bar\phi$ has been chosen so that $x$ lies in a fixed compact set $C\subset X$
(which may be taken to be the fundamental domain for the action of $\rho(\Gamma)$ on the convex hull
of the limit set of $\rho(\Gamma)$). Therefore,  for each $i$, $\rho(\bar\phi(\alpha_i))$ 
translates $x$ by at most $2M+JK+B$, so there are only finitely many choices for
$\bar\phi(\alpha_i)$ for each $i$ and hence only finitely many choices for $\bar\phi$.
This completes the proof in the case that $\Gamma$ acts cocompactly on the symmetric space associated to
a rank one Lie group.

\bigskip

One may mimic this proof in the case that $\Gamma$ is a torsion-free $\delta$-hyperbolic group, 
although several difficulties appear and we will not be able to choose $\mathcal{B}$ to have
the same simple form as above. This is a reasonably straightforward quasification, but as is often
the case, the details and the manipulations of the constants is somewhat intricate.

First, elements of $\Gamma$ need not have axes in $C_\Gamma$. 
However, there exists $(L,A)$ such that every element of $\gamma$ admits a $(L,A)$-quasi-axis $A_\gamma$,
i.e. a $(L,A)$-quasigeodesic which is preserved by $\gamma$ such that if $x\in A_\gamma$, then $d(x,\gamma(x))=||\gamma||$
(see Lee \cite{lee-thesis}).
Second, the nearest point projection need not be well-defined or $1$-Lipschitz.
However, see Gromov  \cite[Lemma 7.3.D]{gromov}, there exists $T$ such that
the nearest point projection $p$ onto any $(L,A)$-quasigeodesic line (or line segment)  in 
$C_\Gamma$ is well-defined up to a distance $T$ and is $(1,T)$-quasi-Lipschitz, i.e.
$$d(p(x),p(y))\le d(x,y)+T$$ 
for all $x,y\in C_\Gamma$.

Given $(L,A)$ and $\delta$, there exists $W>0$ so that if $T$ is a triangle in $C_\Gamma$, whose
edges $e_1$, $e_2$ and $e_3$ are $(L,A)$-quasi-geodesic segments, then there exists a point $z_T\in C_\Gamma$
such that $d(z_T,e_i)<W$ for each $i=1,2,3$.  Moreover, we may assume that $T$ is $W$-slim, i.e. if $x\in e_1$, then
then there exists $y\in e_2\cup e_3$ so that $d(x,y)<W$.
(This follows from Theorem III.H.1.7 and Proposition III.H.1.17
in Bridson-Haefliger \cite{bridson-haefliger}.)

Choose $n>0$ so that if $\gamma\in\Gamma$, then $||\gamma^n||>36(T+2W)$
(which is possible by Theorem III.$\Gamma$.3.17 in Bridson-Haefliger \cite{bridson-haefliger}).
Let $\mathcal B$ be the set of all elements of the form $\alpha_i^n$ or $\alpha_i^n\alpha_j^n$,
where $i\ne j$ and $\alpha_i,\alpha_j\in S$.

Now suppose that $||\phi(\beta)||\le K$ for all $\beta\in{\mathcal{B}}$. We must show
that this determines $\phi$ up to finite ambiguity. Let $\bar\phi$ be a representative
of $\phi$ and, for each $i$,  let $A_i=A_{\bar\phi(\alpha_i^n)}$.
Let $p_i:C_\Gamma\to A_i$ be the nearest-point projection onto $A_i$. We may assume  that
$p_i$ is $\bar\phi(\alpha_i^n)$-equivariant.

For each  distinct pair $i,j\in\{1,\ldots,n\}$,
 let  $x_{ij}\in A_i$ and $y_{ij}\in A_j$ be points such
that $d(x_{ij},y_{ij})=d(A_i,A_j)$.
Let $p_{ij}$ denote nearest point projection onto  a geodesic segment $\overline{x_{ij}y_{ij}}$
joining $x_{ij}$ to $y_{ij}$.
We may assume that  $p_i(\overline{x_{ij}y_{ij}})=x_{ij}$. 

We first claim that 
\begin{equation}
\label{proj-bound1}
p_{ij}(A_i)\subset B(x_{ij},6W).
\end{equation} 
Let $v\in A_i$ and let $v_{ij}=p_{ij}(v)$. Consider the triangle whose sides are $e_1=\overline{v v_{ij}}$, 
the segment $e_2$ of $A_i$ joining $v$ to $x_{ij}$, and $e_3=\overline{x_{ij}y_{ij}}$
(which we may assume is contained in $\overline{x_{ij}y_{ij}}$).
There exists a point $z\in C_\Gamma$ which lies within $W$ of each edge of the triangle.
Let $z_i$ be a point on $e_i$ such that $d(z,e_i)< W$. Then, $d(z_1,z_3)< 2W$. Therefore, since $p_{ij}$ is the
nearest point retraction onto $\overline{x_{ij}y_{ij}}$, $d(z_1,v_{ij})<2W$. So, $d(v_{ij},z_3)<4W$.
On the other hand, since
$d(z_2,z_3)<2W$ and $p_i(\overline{x_{ij}y_{ij}})=x_{ij}$, we see that $d(z_3,x_{ij})<2W$. Therefore,
$d(v_{ij},x_{ij})<6W$ which completes the proof of  (\ref{proj-bound1}).

Let $x_{ij}^\pm$ be two points on $A_i$ such that $d(x_{ij}^\pm,x_{ij})=12(T+2W)$ and if $y$ lies on $A_i$
and is not between $x_{ij}^+$ and $x_{ij}^-$, then $d(z,x_{ij})>12(T+2W)$. We next claim that
if \hbox{$Q_{ij}=p_i^{-1}(\overline{x_{ij}^{+} x_{ij}^{-}})$,} then
\begin{equation}
\label{proj-bound2}
p_{ij}(C_\Gamma-Q_{ij})\subset B(x_{ij},C_0)\ \ {\rm where}\ \ \ C_0=6(T+2W).
\end{equation}
This is the key observation needed to generalize the above proof.
If not, there exists $w$ such that if
$w_i=p_i(w)$ and $w_{ij}=p_{ij}(w)$, then  \hbox{$d(w_i,x_{ij}) > 12(T+2W)$} and $d(w_{ij},x_{ij})>C_0$.
We may assume that \hbox{$p_i(\overline{ww_i})=w_i$} and that  \hbox{$p_{ij}(\overline{ww_{ij}})=w_{ij}$}. Let $e_w$ be
the segment of $A_i$ joining $w_i$ to $x_{ij}$.
Since $p_{ij}$ is
$(1,T)$-quasi-Lipschitz and $p_{ij}(A_i)\subset B(x_{ij},6W)$,
$$d(e_w,\overline{ww_{ij}})\ge C_0-T-6W=5T+6W.$$

Choose a point $z\in e_w$ such that 
$$d(z,x_{ij})>6(T+2W)\ \ \textrm{and}\ \ d(z,w_i)>6(T+2W).$$
Since $d(x_{ij},p_{ij}(z))<6W$ and we may assume that $\overline{x_{ij}w_{ij}}\subset \overline{x_{ij}y_{ij}}$, 
$$d(z,\overline{x_{ij}w_{ij}})>6(T+2W)-6W=6(T+W).$$
First consider the triangle formed by $e_w$, $\overline{x_{ij}w_{ij}}$
and $\overline{w_iw_{ij}}$.
Since the triangle is $W$-slim and $d(z,\overline{x_{ij}w_{ij}})>6(T+W)$, there exists
$z'\in \overline{w_ix_{ij}}$ such that $d(z,z')<W$. Now consider the geodesic triangle formed by 
$\overline{ww_{ij}}$, $\overline{ww_i}$
and $\overline{w_iw_{ij}}$. Since this triangle is also $W$-slim 
and 
$$d(z',\overline{ww_{ij}})>d(e_w,\overline{ww_{ij}})-W\ge 5(T+W),$$
there exists
$z''\in \overline{ww_i}$ such that $d(z',z'')<W$, so $d(z,z'')<2W$. But one easily sees that
this contradicts the facts that $p_i(w)=p_i(z'')=w_i$ and $d(z,w_i)>6(T+2W)$. This completes the proof
of (\ref{proj-bound2}).

Notice that we may choose 
$x_{ji}=y_{ij}$ and $y_{ji}=x_{ij}$. We define \hbox{$R_{ij}=p_j^{-1}(\overline{x_{ji}^+x_{ji}^-})=Q_{ji}$.}
So, (\ref{proj-bound2}) implies that 
$$p_{ij}(C_\Gamma-R_{ij})\subset B(y_{ij},C_0).$$

Therefore,  if 
$$d(A_i,A_j)=d(x_{ij},y_{ij})>C_1=2C_0+K+T,$$
then, again since $p_{ij}$ is $(1,T)$-quasi-Lipschitz,
$C_\Gamma-Q_{ij}$ and $C_\Gamma-R_{ij}$ are disjoint and 
$$d(C_\Gamma-Q_{ij},C_\Gamma-R_{ij})>K.$$
Since $||\bar\phi(\alpha_i^n)||>36(T+W+\delta)$ and $p_i$ is $\bar\phi(\alpha_i^n)$-equivariant,
$Q_{ij}$ is contained in a fundamental domain for the action of $<\bar\phi(\alpha_i^n))>$.
Similarly, $R_{ij}$ is
contained in a fundamental  domain for the action of $<\bar\phi(\alpha_j^n)>$,
$Q_{ij}\cap R_{ij}$ is contained in a fundamental domain for the action of the  subgroup of $\rho(\Gamma)$
generated by $\bar\phi(\alpha_i^n)$ and $\bar\phi(\alpha_j^n)$. Moreover, $Q_{ij}\cap R_{ij}$ is
contained in a fundamental domain for $\bar\phi(\alpha_i^n\alpha_j^n)$ 
and the fixed points of $\bar\phi(\alpha_i^n\alpha_j^n)$
are contained in $\overline{X-Q_{ij}}$ and $\overline{X-R_{ij}}$.
Since $d(C_\Gamma-Q_{ij},C_\Gamma-R_{ij})>K$,
$\|\bar\phi(\alpha_i^n\alpha_j^n)||>K$. Since, we have assumed that $||\bar\phi(\alpha_i^n\alpha_j^n)||\le K$,
we may conclude that $d(A_i,A_j)\le C_1$.

The remainder of the proof proceeds very much as in the convex cocompact case. Notice that the common perpendicular
joining $A_i$ to $A_j$ continues to be replaced with a geodesic joining $x_{ij}$ to $y_{ij}$. 

\medskip\noindent
{\bf Remark:} It would be interesting to know whether or not, in the general case of a torsion-free
word hyperbolic group, one can choose $\mathcal{B}$ to
consist simply of  every element in some finite generating set and the products of any two
elements of the generating set.

\end{document}